\documentclass[11pt]{amsart}
\usepackage{amssymb}
\newtheorem*{mytheorem}{\mytheoremname}
\newcommand{\mytheoremname}{foo}
\newenvironment{myproclaim}[1]{\renewcommand{\mytheoremname}{#1}\begin{mytheorem}}{\end{mytheorem}}
\def\a{\alpha}
\def\b{\beta}
\def\d{\delta}
\def\w{\omega}
\begin{document}

\title{Various topologies on trees}
\author{Peter J. Nyikos}
\address{Department of Mathematics, University of South Carolina, Columbia 
SC 29208}
\email{nyikos@math.sc.edu}
\thanks{Research partially supported by NSF Grant DMS-9322613.}
\thanks{This article was published as {\it Various topologies on trees}, 
in: Proceedings of the Tennessee Topology Conference, P.~R.~Misra and 
M.~Rajagopalan, eds., World Scientific Publishing Co., 1997, pp.~167--198.}
\begin{abstract}
This is a survey article on trees, with a modest number of proofs to give 
a flavor of the way these topologies can be efficiently handled. Trees are 
defined in set-theorist fashion as partially ordered sets in which the 
elements below each element are well-ordered.  A number of different 
topologies on trees are treated, some at considerable length.  Two 
sections deal in some depth with the coarse and fine wedge topologies, and 
the interval topology, respectively.  The coarse wedge topology gives a 
class of supercompact monotone normal topological spaces, and the fine 
wedge topology puts a monotone normal, hereditarily ultraparacompact 
topology on every tree.  The interval topology gives a large variety of 
topological properties, some of which depend upon set-theoretic axioms 
beyond ZFC.  Many of the open problems in this area are given in the last 
section.
\end{abstract}
\maketitle

\section{Trees as abstractions}

Trees, in the everyday sense of the word, generally have a property that 
lends itself almost irresistibly to abstraction.  This is the property of 
repeated branching without rejoining: once the trunk branches off, once a 
branch branches further, etc. there is no subsequent re-combination.

Hence, one hears of such abstractions as ``phone trees,'' ``decision 
trees,'' and ``phylogenetic trees,'' which are based on this, as well as 
``family trees'' which usually conform to it if they just list the near 
ancestors of a person or that person's near descendants. Variants of the 
Latin word ``ramus'', such as ``ramified'' and ``ramification'', are also 
frequently employed in such abstractions.

The abstraction that will take up most of this article is the one used by 
most researchers in set theory and allied fields, with one insignificant 
(but confusing if one is not alert) variant.

\smallskip
\noindent{\bf 1.1. Definition.}  A {\it tree} is a partially
ordered set in which the predecessors of any element are
well-ordered.  [Given two elements $x < y$ of a poset, we
say $x$ is a {\it predecessor} of $y$ and $y$ is a {\it successor} of 
$x$.]
\smallskip

The insignificant variant is that logicians generally put ``successors'' 
in place of ``predecessors'', prompting comments that they are really 
talking about ``root systems''.

It follows from Definition 1.1 that each tree has a set of minimal 
members, above which every member of the tree is to be found.  The 
botanical language continues with:

\smallskip
\noindent{\bf 1.2. Definition.}  If a tree has only one minimal
member, it is said to be {\it rooted} and the minimal member
is called the {\it root} of the tree.  Maximal members (if any)
of a tree are called {\it leaves}, and maximal chains are 
called {\it branches}.
\smallskip

[Recall that a {\it chain} in a poset is a totally ordered subset.  There 
is some conflict in the usage of ``antichain'' where partially ordered 
sets in general are concerned, but fortunately they coincide for trees: an 
{\it antichain} in a tree is a set of pairwise incomparable elements.]

Note the slight deviation from everyday talk: a branch always goes down to 
the bottom level of the tree.  There are standard notations for the levels 
of any tree; the main versions are the one adopted here and the one that 
uses subscripts, putting $T_\a$ where we will use $T(\a)$.

\smallskip
\noindent{\bf 1.3. Notation.}  
If $T$ is a tree, then $T(0)$ is its set of minimal members.  Given an 
ordinal $\a$, if $T(\b)$ has been defined for all $\b < \a$, then $T 
\restriction \a = \bigcup \{T(\b): \b < \a\}$, while $T(\a)$ is the set of 
minimal members of $T \setminus T \restriction \a$.  The set $T(\a)$ is 
called {\it the $\a$-th level of $T$}.
\smallskip

We use the usual notation for intervals, such as $[s, t)  = \{x \in T: s 
\le x < t\}$, and we also adopt the following suggestive notation.

\smallskip
\noindent{\bf 1.4. Notation.}  
Given elements $s < t$ of a tree $T$, let $V_t = \{s \in T: s \geq t\}$, 
and we let $\hat t = \{s \in T: s \le t\}$, given $A \subset T$, let $V_A 
= \bigcup \{ V_a : a \in A\}$ and let $\hat A = \bigcup \{\hat a : a \in 
A\}$.

\smallskip

\noindent{\bf 1.5. Definition.} 
The {\it height} of a tree $T$ is the least ordinal $\a$ such that $T(\a)= 
\emptyset$. Given a cardinal $\kappa$ and an ordinal $\a$, the {\it full 
$\kappa$-ary tree of height $\a$} is the tree of all transfinite 
sequences $f\colon \b \rightarrow \kappa$, for some ordinal $\b < \a$, 
and the order on the tree is end extension: $f \leq g$ iff 
$\operatorname{dom} f \subset \operatorname{dom} g$ and 
$g \restriction \operatorname{dom} f = f$.
\smallskip

The numberings involved in this definition are a bit tricky.  The full 
binary tree of height $\w$ has no elements at level $\w$ and all its 
elements are finite sequences of 0's and 1's.  The full binary tree of 
height $\w +1$, also known as the Cantor tree, has members on its top 
level which are ordinary sequences: there is no $\w$-th term, let alone an 
$\w + 1$-st term.  There are trees of height $\w_1$ with no branches of 
length $\w_1$, such as the tree of ascending sequences of real numbers, 
ordered by end extension.  There are also easy examples of trees of height 
$\omega$ with no infinite branches.

In drawing diagrams of trees, it is traditional to draw line segments 
joining elements to their immediate successors.  These lines are usually 
not meant to be parts of the trees; if they are, point-set topologists 
usually call the resulting objects ``road spaces'' [see Figure 1].  For 
example, what Steen and Seebach [26] refer to as the ``Cantor tree'' is 
more usually called the ``Cantor road space,'' and is formed from the 
Cantor tree by this process of joining successive elements with unit 
intervals. The Moore road space, which Steen and Seebach mention but do 
not define explicitly, can be formed from the Cantor road space by adding 
copies of the unit interval as successors to each point on the $\w$-th 
level of the Cantor tree.

\begin{figure}
\begin{picture}(320,150)(0,0)
\thinlines
\put(160,0){\circle*{2}}
\put(164,0){\scriptsize $\langle \rangle$}

\put(80,80){\circle*{2}}
\put(84,80){\scriptsize $\langle 0\rangle$}
\put(240,80){\circle*{2}}
\put(244,80){\scriptsize $\langle 1\rangle$}

\put(40,120){\circle*{2}}
\put(44,120){\scriptsize $\langle 0,0\rangle$}
\put(120,120){\circle*{2}}
\put(124,120){\scriptsize $\langle 0,1\rangle$}
\put(200,120){\circle*{2}}
\put(204,120){\scriptsize $\langle 1,0\rangle$}
\put(280,120){\circle*{2}}
\put(284,120){\scriptsize $\langle 1,1\rangle$}

\put(20,140){\circle*{2}}
\put(23,140){\scriptsize $\langle 0,0,0\rangle$}
\put(60,140){\circle*{2}}
\put(63,140){\scriptsize $\langle 0,0,1\rangle$}
\put(100,140){\circle*{2}} 
\put(103,140){\scriptsize $\langle 0,1,0\rangle$} 
\put(140,140){\circle*{2}} 
\put(143,140){\scriptsize $\langle 0,1,1\rangle$} 
\put(180,140){\circle*{2}}
\put(183,140){\scriptsize $\langle 1,0,0\rangle$}
\put(220,140){\circle*{2}}
\put(223,140){\scriptsize $\langle 1,0,1\rangle$}
\put(260,140){\circle*{2}}
\put(263,140){\scriptsize $\langle 1,1,0\rangle$}
\put(300,140){\circle*{2}}
\put(303,140){\scriptsize $\langle 1,1,1\rangle$}

\put(10,150){\circle*{2}}
\put(30,150){\circle*{2}}
\put(50,150){\circle*{2}}
\put(70,150){\circle*{2}}
\put(90,150){\circle*{2}} 
\put(110,150){\circle*{2}} 
\put(130,150){\circle*{2}} 
\put(150,150){\circle*{2}} 
\put(170,150){\circle*{2}}
\put(190,150){\circle*{2}}
\put(210,150){\circle*{2}}
\put(230,150){\circle*{2}}
\put(250,150){\circle*{2}}
\put(270,150){\circle*{2}}
\put(290,150){\circle*{2}}
\put(310,150){\circle*{2}}

\put(160,0){\line(1,1){150}}
\put(160,0){\line(-1,1){150}}
\put(80,80){\line(1,1){70}}
\put(240,80){\line(-1,1){70}}
\put(40,120){\line(1,1){30}}
\put(120,120){\line(-1,1){30}}
\put(200,120){\line(1,1){30}}
\put(280,120){\line(-1,1){30}}
\put(20,140){\line(1,1){10}}
\put(60,140){\line(-1,1){10}}
\put(100,140){\line(1,1){10}}
\put(140,140){\line(-1,1){10}}
\put(180,140){\line(1,1){10}}
\put(220,140){\line(-1,1){10}}
\put(260,140){\line(1,1){10}}
\put(300,140){\line(-1,1){10}}
\end{picture}
\caption{The Cantor tree (or the Cantor road space depending on how you 
interpret it).}
\end{figure}
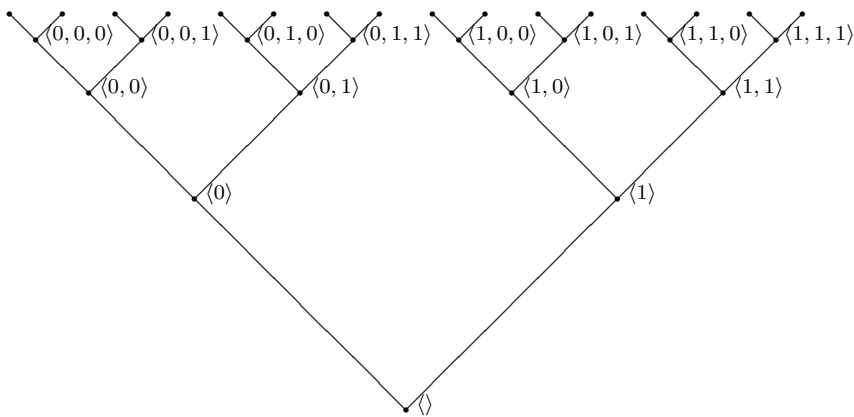

In other branches of topology, and elsewhere in mathematics and the 
sciences, it may be a different story. One popular definition among 
topologists (cf.~[19]) is ``simply connected graph,'' where a graph is 
defined as a nonempty connected 1-complex.  For such spaces, ``simply 
connected'' is equivalent to any two distinct points being the endpoints 
of a unique arc.

Biologists use ``trees'' in a similar way.  Some biologists use the line 
segments in their phylogenetic trees to represent the actual species 
studied, while each fork in their trees represents a speciation event.  
Topologically, they treat their trees as though they were subsets of the 
plane, and they quite correctly observe that the topology depends on the 
actual evolutionary events.  Indeed, once the root of the tree is 
identified, two topologies are equal if and only if they depict the same 
phylogenetic relationship of the species defined.

To minimize confusion between this kind of tree and the trees of 
Definition 1.1, I will only use the word ``tree'' in this other way one 
more time, at the end of the following section.

\section{Comparison of topologies and their convergence properties}

There are many topologies which flow naturally out of the order structure 
of trees.  The ones we discuss have fairly straightforward generalizations 
to partially ordered sets, and some have already been applied to more 
general kinds of posets.

Examples 1 through 4 below will receive additional attention in later 
sections and so we will concentrate on aspects of their convergence in 
this section. We will only need a few concepts pertaining to convergence 
since the topologies we consider are all quite well behaved.

\smallskip
\noindent{\bf 2.1. Definition}  A space $X$ is {\it Fr\'echet-Urysohn}
{\rm [resp. {\it radial\/}]}
if, whenever a point $x$ is in the closure of $A \subset X$, there
is a sequence [{\it resp.} a well-ordered net] in $A$ converging to
$x$.  
\smallskip

A pair of more general concepts will be defined below (Definition 2.3.).

We begin with the topology that is often simply called ``the tree topology''
by set-theoretic topologists.  

\begin{myproclaim}{Example 1}  The {\bf interval topology} on a tree $T$ is the one
whose base is all sets of the form $(s, t] = \{x \in T: s < x \le t\}$, 
together with all singletons $\{t\}$ such that $t$ is a minimal member of $T$.
\end{myproclaim}

It is easy to see that every tree is radial in the interval topology, and 
that a tree is Fr\'echet-Urysohn in this topology iff it is of height $\le 
\w_1$ iff it is first countable.  In fact, a point $t$ is an accumulation 
point of a set $A$ if, and only if, it is in the closure of $A \cap \hat 
t$, and we can order the elements of this set in their natural order to 
produce a well-ordered net converging to $t$, and if it is of countable 
cofinality then we have a sequence converging to $t$.

Every tree is locally compact in the interval topology, and is Hausdorff 
(hence Tychonoff, and 0-dimensional) iff its pseudo-suprema of nonempty 
chains all consist of one point:

\smallskip
\noindent{\bf 2.2 Definition.}  Given a chain $C$ in a tree $T$
such that $C$ is bounded above, the {\it pseudo-supremum of
$C$} is the set of minimal upper bounds of $C$; in other words,
the set of minimal members of $\{t \in T: c \le t \text{ for all }
c \in C \}$.
\smallskip

Pseudo-suprema are always nonempty because they are only defined for 
chains that are bounded above. By the usual conventions, the 
pseudo-supremum of the empty chain is the bottom level $T(0)$ of $T$.  
This is called the trivial pseudo-supremum, and all pseudo-suprema of 
nonempty chains will be called non-trivial.

\begin{myproclaim}{Example 2} The {\bf fine wedge topology} on a tree
 is the topology whose subbase is the collection of all
sets $V_t$ and their complements. 
\end{myproclaim}

It is easy to see that a local base at each point $t$ is formed
by sets of the form
$$
W_t^F = V_t \setminus \bigcup \{V_s : s \in F\} = V_t \setminus V_F
$$
where $F$ is a finite set of successors of $t$.  Of course,
we can restrict ourselves to immediate successors for membership
in $F$. 

It follows from this that a point is isolated in the fine wedge topology 
iff it has at most finitely many immediate successors, and is a point of 
first countability iff it has at most countably many immediate successors.  
But in any case, the topology is always Fr\'echet-Urysohn.  Indeed, $t$ is 
in the closure of $A$ iff $A$ meets $V_s$ for infinitely many immediate 
successors $s$ of $t$, and a sequence $\langle a_n: n \in \w \rangle$ in 
$A$ will converge to $t$ iff only finitely many $a_n$ are above the same 
successor of $t$ and only finitely many are outside $V_t$.

The name for the following topology is inspired by the shape of its
basic open sets.

\begin{myproclaim}{Example 3} The {\bf chevron topology} on a tree $T$ is the 
one whose base
consists of all $\{m\}$ such that $m$ is minimal in $T$, together
with all sets of the form 
$$
C[s, t] = (V_s \setminus V_t) \cup \{t\}
$$ 
such that $s \le t$, where $s$ is either 
minimal or on a successor level.
\end{myproclaim}

Every tree is radial in the chevron topology.  Indeed, a point $t$ is in 
the closure of $A$ if, and only if, either $t \in A$ or $t$ is on a limit 
level and $A$ meets $V_x \setminus V_t$ for cofinally many $x \in \hat t 
\setminus \{t\}$.  In the latter case, we can select, for each $x < t$, an 
element $a_x \in A$ such that $a_x \in V_y \setminus V_t$ for some $y \in 
[x, t)$ and then the well-ordered net $\langle a_x: x <t \rangle$ 
converges to $t$.

It is easy to see that the characters of points are the same in the 
interval and chevron topologies; in particular, the same points are 
isolated in both topologies. Of course, the chevron topology is coarser 
than the interval topology, and strictly coarser in many trees, such as 
Cantor tree, where the top level is easily seen to be homeomorphic to the 
Cantor set in the chevron topology.  In the lattice of all topologies on 
$T$, the least upper bound [called the {\it join}] of the chevron topology 
with the fine wedge topology is the discrete topology since $t$ is the 
only point in $C[s,t] \cap V_t$.

The next four topologies all coincide for trees in which nontrivial 
pseudo-suprema are suprema.  They also agree on the relative topology 
which results when all non-trivial pseudo-suprema of more than one point 
are removed from the tree.  [This should not be confused with the topology 
on the resulting tree that satisfies the formal definition of the 
respective examples.] Example 4b is the coarsest possible topology which 
produces such agreement, while Example 4c is the finest.

\begin{myproclaim}{Example 4a} The {\bf split wedge topology} is the 
greatest lower bound {\rm [i.e., the meet]} of the chevron and fine wedge 
topologies in the lattice of topologies on $T$.
\end{myproclaim}

Note that in a finitary tree [that is, one in which no element has 
infinitely many immediate successors] the split wedge topology and chevron 
topology coincide: the tree is simply discrete in the fine wedge topology.

In other trees, we can construct local bases at each point in the split 
wedge topology by letting their members be simply the union of a basic 
chevron neighborhood and a basic fine wedge neighborhood. Indeed, the 
resulting set is open in both topologies, hence in their meet.  Because of 
this, every tree is radial in the split wedge topology: given $A$ with $t$ 
in its closure, $t$ must be in the closure either of $A \setminus V_t$ or 
of $A \cap V_t$; and then we follow the argument for the corresponding 
finer topology.

The following topology differs from the split wedge topology only in that 
sets of pseudo-suprema [except for the trivial pseudo-supremum $T(0)$] are 
indiscrete rather than discrete in the relative topology.  This allows a 
third possibility for points to be in the closure of $A$, but every 
constant net in an indiscrete space converges to every point in the space, 
so the following topology is again radial.

\begin{myproclaim}{Example 4b} The {\bf coarse wedge topology} on a tree $T$
is the one whose subbase is the set of all wedges $V_t$ and
their complements, where $t$ is either miminal or on a successor
level. 
\end{myproclaim}

If $t$ is minimal or on a successor level, then a local base is formed by 
the sets $W_t^F$ exactly as in the fine wedge topology, with $F$ a finite 
set of immediate successors of $t$.  If, on the other hand, $t$ is on a 
limit level, then one must use $W_{s}^F$ such that $s$ is on a successor 
level below $t$. However, the most appropriate $F$ to take are not sets of 
immediate successors of $s$ but sets of immediate successors of $t$.  
Given any $W_{s}^F$ containing $t$, one can find $t' \in [s, t)$ so that 
the only members of $F$ above $t'$ are those above $t$, and then 
$W_{t'}^G$ is of this form, where $G = F \cap V_t$.

An attractive feature of the coarse wedge topology is that it always has a 
base of clopen sets, even if some nontrivial pseudo-suprema are not 
suprema.  The fine wedge topology and the next example are the only other 
ones that have this feature, of the topologies considered here.

\begin{myproclaim}{Example 4c} The {\bf Lawson topology} on a tree $T$
is the one whose subbase is the set of all wedges $V_t$ and
their complements, where $t$ is not the supremum of a nonempty 
chain in $\hat t \setminus \{t\}$. \end{myproclaim} 

The Lawson topology and the fine wedge topologies are the only ones which 
are invariably Hausdorff for all trees.  The Lawson topology is radial, by 
the same argument as for the split wedge topology. Of course, points in 
pseudo-suprema with more than one element are more easily handled in the 
Lawson topology, because they are isolated.

Another closely related topology is intermediate between the split wedge 
and coarse wedge topologies, giving sets of pseudo-suprema the cofinite 
topology.  Since every injective $\w$-sequence converges to each point of 
a space with cofinite topology, this topology too is radial:

\begin{myproclaim}{Example 4d} The {\bf hybrid wedge topology} on a tree
$T$ is the one whose subbase consists of all complements of 
wedges $V_t$ together with those wedges $V_s$ for which $s$
is either minimal or a successor. \end{myproclaim}

Note that doing it the other way around---all wedges plus complements
of wedges based on successors or mimimal members---simply produces
the fine wedge topology because the basic sets $W_t^F$ such that $F$
consists of immediate successors of $t$, are all there.

So far, the topologies we have been considering are all Hausdorff and 
zero-dimensional if all nontrivial pseudo-suprema are singletons, hence 
suprema (and we can drop the conditional clause for the fine wedge and 
Lawson topologies).  With one exception (Example 7) this is not the case 
with the remaining topologies of this section. These remaining topologies 
will not be dealt with in subsequent sections and the reader may skip to 
Section 3 now or later with no loss of continuity.

The next two topologies are Hausdorff iff they are $T_1$ iff no element is 
above any other element, i.e. if $T(0)$ is all of $T$; and in this case, 
they are discrete.  The first one can be thought of as `one half of the 
Lawson topology':

\begin{myproclaim}{Example 5} The {\bf Scott topology} on a tree $T$ is
the one for which sets of the form $V_t$ are a base, where $t$ is not the
supremum of a nonempty chain in $\hat t \setminus \{t\}$. \end{myproclaim}

For arbitrary posets, one has to use a different description, easily seen 
equivalent for trees: the Scott-open subsets of a poset $P$ are those 
upper sets $U$ such that no member is the supremum of a directed subset of 
$P \setminus U$. [A subset $U$ of a poset is said to be an {\it upper set} 
if $V_p \subset U $ whenever $p \in U$.  Thus, for example, the upper 
subsets of $\mathbb{R}$ are the right rays, and those of the form $(a, 
+\infty)$ are the Scott open subsets.]

Of course, every tree is a $T_0$-space in the Scott topology. It is a 
radial space by the same arguments that apply to the Lawson topology, only 
they are simpler since any constant sequence in $V_t$ converges to $t$ in 
the Scott topology. This applies {\it a fortiori} to the next topology, 
where a point $t$ is in the closure of $A$ iff $A$ meets $V_t$.

\begin{myproclaim}{Example 6} The {\bf Alexandroff discrete topology} is 
the one for which all sets of the form $V_t$ form a base.
\end{myproclaim}

Examples 5 and 6 have the propery that the order can be recovered from the 
topology by setting $x \le y$ iff $x$ is in the closure of $\{y\}$.

While these last two topologies may be ``uninteresting'' from the point of 
view of most general topologists, they have great significance from other 
points of view.  The Alexandroff discrete topology, generalized to posets, 
is the one behind words like ``open'' and ``dense'' in the applications of 
forcing.

Forcing is a method of producing models of set theory, pioneered by Paul 
Cohen, who used it in 1963 to show that the continuum hypothesis is 
independent of the usual axioms of set theory.  It has revolutionized set 
theory and a number of other branches of mathematics, especially 
set-theoretic topology and the theory of Boolean algebras.

The Scott topology is important in theoretical computer science 
(cf.~[18]).  Appropriately enough, it is named after the leading computer 
scientist Dana Scott, who showed [25] that continuous lattices equipped 
with this topology are precisely the injective objects in the category of 
$T_0$-spaces and continuous functions.

When the discrete topology and one more topology are added, and we 
restrict our attention to trees in which pseudo-suprema are suprema, the 
foregoing topologies form a sublattice of the lattice of all topologies, 
as shown in Fig. 2. The pentagon on the right shows that this is not a 
modular lattice.

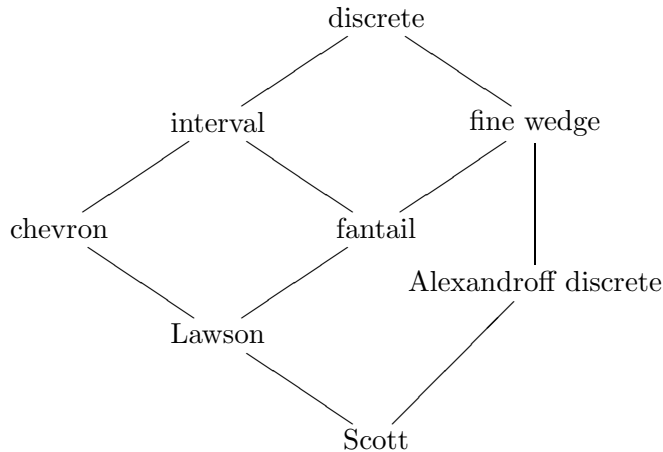
\begin{figure}[Lattice of topologies]
\begin{picture}(220,160)(0,0)
\thinlines
\put(120,160){\makebox(0,0){discrete}}
\put(60,120){\makebox(0,0){interval}}
\put(180,120){\makebox(0,0){fine wedge}}
\put(0,80){\makebox(0,0){chevron}}
\put(120,80){\makebox(0,0){fantail}}
\put(180,60){\makebox(0,0){\centering Alexandroff discrete}}
\put(60,40){\makebox(0,0){Lawson}}
\put(120,0){\makebox(0,0){Scott}}
\put(171,126){\line(-3,2){40}}
\put(69,126){\line(3,2){40}}
\put(129,86){\line(3,2){40}}
\put(111,86){\line(-3,2){40}}
\put(9,86){\line(3,2){40}}
\put(180,66){\line(0,1){46}}
\put(69,46){\line(3,2){40}}
\put(51,46){\line(-3,2){40}}
\put(126,6){\line(1,1){46}}
\put(111,6){\line(-3,2){40}}
\end{picture}
\caption{{\bf Lattice of topologies.} This diagram is valid for trees in which 
every nontrivial pseudo-supremum is a supremum. The situation is more 
complicated for general trees: the Scott Topology (and hence also the 
Lawson topology) is not always coarser than the interval topology (nor, 
{\it a fortiori}, than the other topologies below the interval topology in 
the diagram). The meet of the fantail and chevron topologies in general is 
the split wedge topology.}
\end{figure}

I have given the name ``fantail topology'' to the topology that is the 
meet of the interval and fine wedge topologies, because of the pictures I 
associate with the basic neighborhoods as defined below.

\begin{myproclaim}{Example 7} The {\bf fantail topology} is the 
one for which a base is the collection of all sets of the form
$
\bigcup \{W_x^{F(x)}: s \le x \le t\} 
$
such that $s$ is either minimal or a successor,
and $F(x)$ is a finite subset of $V_x$ for all $x$, satisfying the
following condition: if $x < t$ then $x' \in F(x)$, where $x'$ is the 
immediate successor of $x$
which satisifes $x' \le t$.
\end{myproclaim}

This is a more complicated topology than the ones considered so far, and 
it is the only one which is not always radial.  However, it is the next 
best thing in a sense:

\smallskip
\noindent{\bf 2.3. Definition} A space $X$ is {\it pseudo-radial} if
closures can be taken by iterating the operation of taking 
limits of convergent well-ordered nets; the number of iterations
required is the {\it chain-net order of $X$}.
\smallskip

More formally: given a set $A \subset X$, let $A^{\tilde{~} 1}$ be the 
set of all points which are limits of well-ordered nets from $A$. If $\a$ 
is an ordinal and $A^{\tilde{~} \a}$ has been defined, let 
$A^{\tilde{~} {\a +1}} = (A^{\tilde{~} \a})^{\tilde{~} 1}$, 
while if $\a$ is a limit ordinal we define $A^{\tilde{~} \a}$ to be the 
union of all the $A^{\tilde{~} \b}$ such that $\b < \a$. The first 
ordinal $\a$ such that $A^{\tilde{~} {\a +1}} = A^{\tilde{~} \a}$ for 
all $A \subset X$ is the chain-net order of $X$, provided 
$A^{\tilde{~} \a +1} = A^{\tilde{~} \a}$ is always the closure of $A$; 
this condition characterizes pseudo-radial spaces.

\begin{myproclaim}{2.4. Theorem} Every tree is pseudo-radial, of chain
net order $\le 2$ in the fantail topology.  The order is exactly
2 in any tree of height $> \w$ in which every element of $T(n)$ has 
infinitely many immediate successors for $n \in \w$.  \end{myproclaim}

Before showing this, it is helpful to make some observations and to define 
another concept pertaining to general spaces.

We were `fortunate' to have the union of a basic chevron neighborhood
and a basic fine wedge neighborhood be open in both topologies.
In the case of more general spaces, and in particular the case
where the interval topology is substituted for the fantail topology,
we can expect only that such unions form a weak base:

\smallskip
\noindent{\bf 2.5. Definition.}  
Let $X$ be a set.  A {\it weak base on $X$} is a family of filterbases 
$\mathbb{B} = \{\mathcal{B}(x) : x \in X\}$ such that $x \in B$ for all $B 
\in \mathcal{B}(x)$.  The topology induced on $X$ by $\mathbb{B}$ is the 
one in which a set $U$ is open if, and only if, there exists for each 
point $x \in U$ a member $B$ of $\mathcal{B}(x)$ such that $B \subset U$.
\smallskip

Of course, every system of ordinary neighborhood bases is a weak base, and 
this is something worth keeping in mind when reading the following lemma.

\begin{myproclaim}{\bf 2.6. Lemma} 
If $\tau_1$ and $\tau_2$ are topologies on $X$ and $\mathbb{B}_1$ and 
$\mathbb{B}_2$  are  weak bases for $\tau_1$ and $\tau_2$ respectively, 
then the weak base $\mathbb{B}$ in which $\mathcal{B}(x) = \{B_1 \cup B_2: 
B_i \in \mathcal{B}_i(x) \text{ for }
i =  1, 2\}$ is a weak base for the meet of the $\tau_i$.
\end{myproclaim}

\begin{proof}
Let $U$ be open in both $\tau_i$ --- equivalently, in
their meet.  For each $x \in U$ and each $i$ there exists $B_i \in 
\mathcal{B}_i(x)$ such that $B_i(x) \subset U$, so $B_1(x) \cup B_2(x) 
\subset U$.
Conversely, suppose that $V \subset X$ and for each $x \in V$,
there exist $B_i \in \mathcal{B}_i(x) $ such that $B_1 \cup B_2 
\subset V$.
Then $V$ is open in $\tau_i$ for both $i$, hence in their meet.
\end{proof}

Now, in the case where the $\tau_i$ are the fine wedge and interval 
topologies, letting $\mathcal{B}_1(x)$ be the local base of all sets 
$W_x^F$, and letting $\mathcal{B}_2(x)$ be the set of all intervals $(s, 
x]$, gives us a weak base for the fantail topology in which $B_1 \cup B_2$ 
is not always open in the fine wedge topology,
 hence it is not always open in the fantail topology.  However,
by attaching a fine wedge neighborhood to each point of $ B_2 =(s, t]$
we do produce a set that is open in both topologies, and it is easy
to see that every set that is open in both topologies must contain
a set of this form; among these are the basic open sets described in
the statement of Example 7.  They are also easily seen to be closed
in both topologies.

\begin{proof}[Proof of Theorem 2.4]
First we show that every point in the closure of $A$ is the limit of a 
convergent well-ordered net in $\tilde A$ for all $A \subset T$, where 
$\tilde A$ is the set of all limits of convergent sequences in $A$.  This 
we do by showing that if $t$ is in the closure of $A$, then either $t \in 
\tilde A$ or a cofinal subset of $\hat t \setminus \{t\}$ is in $\tilde 
A$.  Of course, there will be a well-ordered net in this cofinal subset 
converging to $t$.

If $t$ is in the closure of $A \cap V_t$, then there is a sequence in $A 
\cap V_t$ converging to $t$ just as in the fine wedge topology.  So 
suppose not; then $t$ is in the closure of $A \setminus V_t;$ this of 
course implies that $t$ is on a limit level.  If there were no cofinal 
subset of $\hat t \setminus \{t\}$ in $\tilde A$, then we could attach 
fine wedge neighborhoods missing $A$ to all members of a final segment 
$[s, t)$], thereby keeping $t$ out of the closure of $A \setminus V_t$; 
and thus we can build a basic neighborhood as in the initial presentation 
of Example 7, missing $A$, contradicting the claim that $t$ is in the 
closure of $A$.

To show that chain-net order is exactly 2 in trees as described in the 
second sentence, we will produce a copy of the Arens space $S_2$ in any 
such tree.

The Arens space can be defined as the space whose underlying set $S_2$ is 
$\{x_\w\} \cup \{x_n : n \in \w\} \cup \{x_k^n : n, k \in \w\}$ faithfully 
indexed, e.g. $ x_n^k = x_m^i$ iff $n = m$ and $k = i$; and in which a 
weak base is given by $\{\mathcal{B}(p) : p \in S_2\}$ where 
$\mathcal{B}(x_k^n) = \{\{x_k^n\}\}$, $\mathcal{B}(x_n)$ is the collection 
of all sets $ \{x_n\} \cup \{x_k^n : k \ge j\}$ as $j$ ranges over $\w$; 
and $\mathcal{B}(x_\w)$ is the collection of all sets $\{x_\w\} \cup \{x_j 
: j \ge n\}$ as $n$ ranges over $\w$.  Some elementary properties of $S_2$ 
are that $\langle x_n : n \in \w \rangle$ converges to $x_\w$; that all 
the points $x_k^n$ are isolated; and that the points $x_n (n \in \w)$ are 
points of first countability with each set $ \{x_n\} \cup \{x_k^n : k \ge 
j\}$ a clopen copy of $\w + 1$.

The most crucial feature, for our purposes, is that the set of isolated 
points has $x_\w$ in its closure, yet no sequence of isolated points can 
converge to $x_\w$.  To see this, note that any sequence that meets some 
`column' $\{x_k^n : k \in \w\}$ in an infinite set cannot converge to 
$x_\w$, and also that any sequence of isolated points that meets each 
`column' in a finite set does not have $x_\w$ in its closure, because we 
can take a member $B$ of $\mathcal{B}(x_\w)$ and attach a clopen set $ 
\{x_n\} \cup \{x_k^n : k \ge j_n\}$ missing the sequence to each $x_n \in 
B$, producing an open neighborhood of $x_\w$ that misses the sequence.

Now it is routine to build a copy of $S_2$ in any tree as described.  Let 
$\{x_\a : \le \w\}$ be represented by $\hat t$ where $t$ is any point on 
level $\w$; of course, $x_\w$ is represented by $t$ itself.  For each $x_n 
(n \in \w)$ the sequence $ \{x_n\} \cup \{x_k^n : k \in \w \}$ is 
represented by some countably infinite set of immediate successors of 
$x_n$ not in $\hat t$.  The weak base given in the paragraph following the 
proof of Lemma 2.6 traces a weak base on this subspace exactly as in 
$S_2$.  Moreover, the resulting copy of $S_2$ has no more points of the 
tree in its closure except perhaps points in the pseudo-supremum of $\{x_n 
: n \in \w\}$; and these do not alter the fact that $t$ cannot have any 
sequence converge to it from the copy of $S_2$ other than those in which a 
cofinite subsequence is in $\{x_n : n \in \w\}$.
\end{proof}

This concludes our treatment of the fantail topology.  The next three 
examples are taken from [18]. Our list of tree topologies in this section 
will be concluded with their meets and joins with each other and with the 
topologies given so far. Two of the topologies are an instance of a 
general motif: given a topology $\tau$ defined on a class of posets, one 
obtains the reverse topology $\tau^d$ on a poset by turning the poset 
upside down (i.e. reversing the order relation), defining $\tau$ on the 
resulting poset, and then turning it back right side up again.

\begin{myproclaim}{Example 8} 
The {\bf Scott${}^{\bf d}$ topology} on a tree $T$ is the one for which 
the subsets $\hat t$ form a base for the topology.
\end{myproclaim}

Indeed, if one inverts a tree, no element is the directed supremum of 
now-lower elements, and so every now-upper set is open.  The 
Scott${}^d$-topology is obviously coarser than the discrete and interval 
topologies and is incomparable to all the remaining topologies.  Its join 
with the Alexandroff discrete topology (and hence with all finer 
topologies) is obviously the discrete topology, and its join with the 
Scott topology (and hence all others not above the Alexandroff discrete 
topology) is clearly the interval topology.

The Scott${}^d$-topology is obviously first countable; in fact, $\{\hat 
t\}$ is a one-member local base at $t$.

\begin{myproclaim}{Example 9}  The {\bf weak topology} on a tree $T$ is
the one in which the sets $\hat t$ form a subbase for the
closed sets. \end{myproclaim}

Since the sets $\hat t$ are downwards closed and linearly ordered, it 
easily follows that any closed set other than the whole space is a finite 
union of sets of the form $\hat t$. Moreover, a point $t$ is in the 
closure of a finite set $F$ iff it is in $\hat F$ iff it is in $\hat s$ 
for some $s \in F$ iff $V_t$ meets $F$. A subset $A$ is dense in $T$ iff 
it either contains an unbounded chain, or the set of suprema of chains in 
$A$ contains an infinite antichain. Hence this topology is radial, and is 
clearly coarser than the Scott topology and hence coarser than all the 
topologies considered so far except the coarse wedge and Scott${}^d$ 
topologies. And its dual is weaker than the Scott${}^d$ topology, of 
course:

\begin{myproclaim}{Example 10.}  
The {\bf weak${}^{\bf d}$ topology} on a tree $T$ is the one in which the 
sets $V_t$ form a subbase for the closed sets. 
\end{myproclaim}

In every tree in which pseudo-suprema are suprema, one can just as easily 
use only those $V_t$ in which $t$ is a minimal or successor element.  In 
any event, it is routine to show that a set is closed in the weak${}^d$ 
topology iff it equals $V_A$ for an antichain $A$ such that $\hat A$ is a 
finitary tree with finitely many minimal members.

This topology is Fr\'echet-Urysohn, with a local base at $t$ consisting of 
sets $T \setminus V_S$ where $S$ is a finite set of points of $T$ that are 
outside of $\hat t$, but have all their predecessors in $\hat t$. This 
includes the points of $T(0) \setminus \hat t$ by the usual conventions on 
the empty set.

The join of the weak${}^d$ topology and the Scott topology is the Lawson 
topology; in fact this is the way the Lawson topology is defined by Lawson 
in [18].  The join of the weak and weak${}^d$ topology is strictly coarser 
than the Lawson topology in general, and it is not hard to see that any 
closed set is a union of a weak${}^d$-closed set [see description above] 
together with finitely many intervals of the form $[s,t]$. 

The meet of the Alexandroff discrete and Scott${}^d$ topologies is the 
indiscrete topology on every rooted tree.  Of course, this also applies 
when a coarser topology replaces either or both of these, but there are 
differences in other trees.  In case of the Alexandroff discrete and 
Scott${}^d$ topologies, we simply have the topological direct sum of the 
indiscrete rooted trees involved.  This is also true if the weak${}^d$ 
topology replaces the Scott${}^d$ topology and/or the Scott topology 
replaces the Alexandroff discrete topology.  On the other hand, the weak 
topology gives the cofinite topology on those trees which consist of the 
single level $T(0)$. Of course, only Examples 8 and 10 are not finer than 
the weak topology, and it is a simple matter to see that in both cases, a 
set is closed in the meet topology if and only if it contains each $V_t$ 
that it meets and is either the entire tree, or else it meets $V_t$ for 
only finitely many minimal $t$.

The remaining meets are more interesting.  They are found in the third 
quadrant of Figure 2, so to speak.

\smallskip
\noindent{\bf Example 11.} 
The meet of
the Scott${}^d$-topology
and the chevron topology can be understood via Lemma 2.6 in the same
way that the fantail topology is.  A weak base at each point $t$ on a 
limit level consists of sets of the form $C[s,t] \cup \hat t$ where 
$C[s, t]$ is a chevron, while
a weak base at a point $s$ on a successor or minimum level consists simply
of  $\{\hat s\}$.  Construction of a base can be worked out as for
the fantail topology.  I picture a typical member informally as a feather
with a wedge cut in the top and finitely many indentations going
all the way to the central shaft, each indentation going a finite
number of steps up the central shaft.  More formally: a local
base at $t$ consists of sets of the form 
$
\hat t \cup (V_{t_0} \setminus \bigcup \{ V_x : x \in A \}) 
$ 
where $t_0$ is the minimal
element of $\hat t$, and $A$ is a finite union of levels of $T$
including the one on which $t$ itself is to be found.
One cannot exclude infinitely many such levels without causing
trouble at the next limit level.  This fact makes it routine to 
show that any tree is Fr\'echet-Urysohn in this topology.

\smallskip
\noindent{\bf Example 12.} 
Similarly, the meet of the 
Scott${}^d$-topology with the fine wedge topology
(and hence with the fantail topology)
can be characterized as the one whose weak base at $t$ is formed 
by attaching sets of the form $W_t^F$ to $\hat t$.  In forming the base,
one could just follow the description of the fantail topology, just making
sure that the basic neighborhoods start with the minimum point $t_0$
of $\hat t$.
Like the fantail topology, this one is pseudo-radial of order $\le 2$.

\smallskip
\noindent{\bf Example 13.} 
Examples 11 and 12 are incomparable; 
their meet has a base formed by taking
a basic Example 11 neighborhood and attaching a set $W_x^F$ to each of
the finitely many points of $\hat t$ in the levels met by $A$.
Of course this is also the meet of the  Scott${}^d$-topology and the
split wedge topology, so it is finer than the weak${}^d$-topology. 

Meets involving the coarse wedge, hybrid wedge, and Lawson topologies are 
left as an exercise for the interested reader.

We close this section with some comments about the Scott topology on 
phylogenetic trees.  This seems to be the appropriate topology for the 
branch of systematics called cladistics, which is centered on those groups 
of organisms which form clades.  Clades are simply sets of organisms 
represented by the various $V_t$ in a phylogenetic tree, and many cladists 
will refuse to even consider taxa that are not clades as legitimate 
scientific entities. Their rationale [which I consider to be inadequate] 
is that one can recover the entire order on the tree by just knowing what 
the clades are.  This is, of course, a very useful thing to be able to do, 
and is very closely related to the fact that the Scott topology allows us 
to recover the order on the tree.

\section{Completeness and Compactness}

In this section we consider some elementary aspects of completeness of a 
tree which depend only on the order, and the kinds of compactness they 
give rise to.

Where trees are concerned, the very fundamental concept of Dedekind 
completeness simply translates to every pseudo-supremum being a supremum. 
Hence it is easy to produce a Dedekind completion for any tree: just give 
every set of pseudo-suprema of more than one element an immediate 
predecessor. Note, however, that the inclusion map of the original tree in 
its Dedekind completion is not an embedding in many of the topologies of 
Section 2 [in particular, not in Examples 1 through 4b nor for 4d; the 
Lawson topology is a noteworthy exception]. This is because the points of 
the original pseudo-suprema are no longer in the closure of the points on 
the earlier levels in most of the topologies. On the other hand, the map 
does have dense range in most of the topologies, including the Lawson 
topology.  Of course, it is always an order-embedding.

Being Dedekind complete is equivalent to the tree being Hausdorff (also to 
being Tychonoff) in the coarse wedge topology and all finer topologies 
considered in Section 2.  It is also equivalent to the space having a base 
of clopen sets in the Lawson, chevron, hybrid, and interval topologies, 
inasmuch as these are always $T_1$.  Since this is not the sort of 
property one usually associates with Dedekind completeness, I am being 
quite sparing of the term in this article where trees are concerned.  
However, in one case it does seem appropriate:

\begin{myproclaim}{3.1 Theorem} 
Let $T$ be a tree.  The following are equivalent.
\begin{itemize}
\item[(i)] $T$ is rooted and Dedekind complete.
\item[(ii)] $T$ is a semilattice downwards; that is,  any two elements
have a greatest lower bound.
\item[(iii)] Every nonempty subset of $T$ has a greatest lower bound.
\end{itemize}
\end{myproclaim} 

\begin{proof}
{\it (i) implies (iii):}  Let $A$ be a nonempty subset
of $T$ and let $B$ be the set of lower bounds for $A$; $B$ is 
nonempty since the tree is rooted and it is clearly a chain.  Since the
levels of the tree are well-ordered, $B$ has a pseudo-supremum on the
first level above which there are no members of $B$, and since it is
a supremum of $B$ it is also the greatest lower bound of $A$.

{\it (iii) implies (ii)}: Obvious.

{\it (ii) implies (i)}: It is obvious that $T$ cannot have more than
one minimal element.  If $T$ had a pseudo-supremum that is not a 
supremum, then any two distinct elements of this pseudo-supremum
would fail to have a greatest lower bound. 
\end{proof}

Actually, the equivalence of (ii) and (iii) is part of a more general 
phenomenon: every infimum is the infimum of some two-element subset.  
Formally:

\begin{myproclaim}{3.2 Theorem} Let $A$ be a set of two or more elements
of a tree $T$ such that $A$ has a greatest lower bound.  Then there are 
elements $a_1$ and $a_2$ of $A$ such that the g.l.b.\ of $a_1$ and $a_2$ 
is the g.l.b.\ of $A$.
\end{myproclaim}

\begin{proof}
Let $t$ be the g.l.b.\ of $A$.  There are at least two distinct
immediate successors of $t$ with elements of $A$ above them, and we choose
$a_1$ and $a_2$ from above two of them. 
\end{proof}

Even where there is no greatest lower bound, one can speak of 
``pseudo-infima'' in analogy with pseudo-suprema.  Then every nonempty 
subset of every tree has a pseudo-infimum, and if $A$ has at least two 
elements, we can find two whose pseudo-infimum is the pseudo-infimum of 
$A$.

The use of ``complete'' does seem quite appropriate in the following 
concepts, and leads to some nice compactness results.

\smallskip
\noindent{\bf 3.3 Definition}  
A tree is {\it branch-complete}
if every branch has a greatest element.  A tree is {\it chain-complete}
if every chain has a supremum.
\smallskip

Branch-completions can trivially be produced by giving any branch a 
greatest element if it does not already have one.  The original tree is 
densely embedded in the resulting tree in Examples 1 through 4d, except 
for the fine wedge topology.  Chain-completions can simply be produced by 
taking a Dedekind completion of a branch completion, or vice versa; if the 
descriptions given above are followed, the same tree results no matter 
which completion is taken first.

The following theorem identifies a rich source of well-behaved
examples of compact Hausdorff spaces.

\begin{myproclaim}{3.4. Theorem} Let $T$ be a tree.
The following are equivalent.
\begin{itemize}
\item[(i)]  $T$ is branch-complete and has finitely many minimal elements.
\item[(ii)] $T$ is compact in the coarse wedge topology.
\end{itemize}
\end{myproclaim}

\begin{proof}
{\it (i) implies (ii):} We will show, in fact,  that if
$T$ is rooted, it
is {\it supercompact}; that is, it has a subbase such that every
open cover has a subcover by two or fewer members.  This implies
compactness by Alexander's subbase theorem.  The result for non-rooted
trees follows immediately since they are topological direct sums
of rooted ones in the coarse wedge topology.

Let $\{V_x : x \in A\}$ and $\{T \setminus V_x : x \in B\}$ be a subbasic 
open cover of $T$.  If $A \cap B \neq \emptyset$, then we immediately have 
a two-member subcover.  Otherwise, pick a member of the cover containing 
the root $t_0$ of the tree. If this is of the form $V_x$ we are done since 
$x = t_0$ and $V_x$ is all of $T$. Otherwise, every member of the cover 
containing $t_0$ is of the form $T \setminus V_x$. If $B$ has a pair of 
incomparable elements, say $x$ and $y$, then $\{T \setminus V_x,\, T 
\setminus V_y \}$ is a subcover.

It remains to consider the case where $B$ is a chain. By 
branch-completeness, $B$ has a pseudo-supremum $P$.  The points of $P$ can 
only be covered by a set of the form $V_a$ with $a < p$ for all $p \in P$. 
But then there exists $b \in B$ such that $b > a$, and then $\{V_a, T 
\setminus V_b\}$ is as desired.

{\it (ii) implies (i)} Since the minimal level of a tree is closed 
discrete in the coarse wedge topology, $T$ can have only finitely many 
elements in this level if it is to be even countably compact.  Also every 
branch must have a maximum member, for if $B$ violates this, $\{T 
\setminus V_b : b \in B, b \text{ is a successor or minimal} \}$ is an 
open cover with no finite subcover.
\end{proof}

\begin{myproclaim}{3.5. Corollary} A tree is compact Hausdorff in the
coarse wedge topology if, and only if, it is chain-complete
and has finitely many minimal elements.  
\hfill $\square$
\end{myproclaim}

The proof of the following theorem will appear in a forthcoming paper.  
Note the absence of any completeness condition.

\begin{myproclaim}{3.6. Theorem}  
Let $T$ be a tree.  The following are equivalent.
\begin{enumerate}
\item $T$ has countably many minimal elements.
\item $T$ is $\w_1$-compact in the coarse wedge topology; that is,
every closed discrete subspace is countable.  
\hfill $\square$
\end{enumerate}
\end{myproclaim}

The foregoing results remain true if the hybrid wedge topology replaces 
the coarse wedge topology, except that the simple proof of 
supercompactness in Theorem 3.4 may fail even if the tree is rooted.  The 
proof of the second implication goes through with no change except for the 
name of the topology.  For the first implication, we make a minor 
modification if $P$ has finitely many elements; in that case, there is the 
additional possibility that finitely many $V_p$ round out the subcover.

The split wedge topology does not add any new compact examples since it 
coincides with the hybrid wedge topology when pseudo-suprema are finite, 
and has an infinite closed discrete subspace otherwise. Of course, this 
applies also to the Lawson topology.  A similar statement holds for the 
chevron topology: it is compact iff it coincides with the hybrid wedge 
topology and the latter is compact.  Equivalently, the tree is finitary 
and every chain has a finite pseudo-supremum.

Theorem 3.4 and Corollary 3.5 have straightforward analogues for countably 
compact spaces.  Proofs are left as an exercise for the reader:

\begin{myproclaim}{3.7. Theorem} Let $T$ be a tree with the coarse wedge
or hybrid wedge topology.
The following are equivalent.
\begin{itemize}
\item[(i)]  $T$ has finitely many minimal elements, and every branch of 
countable cofinality has a maximal element.
\item[(ii)] $T$ is countably compact.
\item[(iii)] $T$ is sequentially compact.
\hfill $\square$
\end{itemize}
\end{myproclaim}

A nice application of the coarse wedge topology was found by Gary 
Gruenhage [11]:

\begin{myproclaim}{3.8 Example} 
A locally compact, metalindel\"of space which is not weakly 
$\theta$-refinable. 
\end{myproclaim}

Let $S$ be a stationary, co-stationary subset of $\w_1$ and let $T$ be the 
set of all compact subsets of $S$, with the end extension order.  Let $X$ 
be the chain-completion of $T$, with the coarse wedge topology.  Then $X^2 
\setminus \Delta$ is as described.  This was the first ZFC example of a 
metalindel\"of regular space that is not weakly $\theta$-refinable.

The following example has found use in functional analysis.  It is 
attributed to D.~H.~Fremlin by Richard Haydon [private communication] and 
has been rediscovered by several researchers, including J.~Bourgain, to 
whom it is attirbuted by J.~Diestel [8, p.~239].

\smallskip
\noindent{\bf 3.9. Example.} Let $\w^*$ stand for the Stone-\v{C}ech
remainder of $\w$; in other words, $\w^* = \b\w - \w$. 
Let $\mathcal{C}_0 = \{\w^*\}$.
Let $\mathcal{C}_1$
be an uncountable collection of disjoint
clopen subsets of $\w^*$ whose union is dense.  
If $\a$ is a successor ordinal and the disjoint 
collection $\mathcal{C}_\a$ of clopen sets
has been defined then $\mathcal{C}_{\a+1}$
is obtained by taking the union of uncountable families of disjoint 
clopen sets in each member of $\mathcal{C}_\a$, each family having 
dense 
union in its respective member.  If $\a$ is a limit ordinal and
$\mathcal{C}_\b$ has been defined
for all $\b < \a$, let $\mathcal{C}_\a$ be the collection of all 
intersections
of maximal chains in $\bigcup \{\mathcal{C}_\b: \b < \a \}$ and let 
$\mathcal{C}_{\a+1}$ be the union of (uncountable) families of 
disjoint
clopen sets in each member of $\mathcal{C}_\a$ whose interior is 
nonempty, 
each family having dense union in the interior of its respective member.  
Continue until a limit ordinal $\gamma$ has been reached such that
every member of $\mathcal{C}_\gamma$ has empty interior, and let 
$\mathcal{T}$ be the tree 
$
\bigcup \{\mathcal{C}_\a: \a \le \gamma \} 
$
ordered by reverse inclusion.  

Of course, $\mathcal{T}$ is chain-complete and rooted, hence compact 
Hausdorff in the coarse wedge topology. What is especially significant is 
that it is homeomorphic in a natural way to the decomposition space of 
$\w^*$ whose elements are the closed nowhere dense sets $F \setminus 
\bigcup \{C \in \mathcal{C}_{\a+1}: C \subset F \} $ as $F$ ranges over 
$\mathcal{C}_\a$ and $\a$ ranges over the ordinals $ \le \gamma$. [Of 
course, if $F \in \gamma$ then $F$ is nowhere dense and 
$\mathcal{C}_{\a+1} = \emptyset$.] The map associating $F \in \mathcal{T}$ 
with this nowhere dense subset is a homeomorphism.  Moreover, if 
$\mathcal{T}$ is a $\pi$-base, then the decomposition map from $\w^*$ to 
the decomposition space is irreducible.  See [2] for the constuction of 
tree $\pi$-bases for $\w^*$ and their uses.

The fact that each member of the decomposition space is nowhere dense 
tells us that no sequence from $\w$ will converge anywhere in the compact 
Hausdorff space which is the quotient space of $\b\w$ formed by 
identifying the members of the decomposition space to points.  So this 
space shares some of the `pathology' of $\b\w$ and yet the set of 
nonisolated points is far better behaved, being radial and having lots of 
convergent sequences.  A few other `nice' properties of the remainder will 
become evident at the beginning of Section 4.  These spaces have been 
studied in an effort to characterize the smallest (uncountable) cardinal 
$\kappa$ such that there is a compact Hausdorff space of cardinality 
$\kappa$ which is compact but not sequentially compact.

A whole class of even `nicer' compactifications is associated in a natural 
way to non-Archimedean spaces:

\smallskip
\noindent{\bf 3.10. Definition.}  A collection $\mathcal{B}$ of
subsets of a set is {\it of rank 1\/} if, given any two members
$B_1$, $B_2$, either $B_1 \cap B_2 = \emptyset$ or $B_1 \subset B_2$
or $B_2 \subset B_1$.  A {\it non-Archimedean space} is a $T_0$
[equivalently, Tychonoff] space with a rank 1 base.
\smallskip

A crucial fact about non-Archimedean spaces is that they actually have a 
base which is a tree under reverse containment ([21]).  This makes the 
proof of such `nice' attributes as ultraparacompactness and 
suborderability very easy, and also leads in a natural way to embedding 
them in compact Hausdorff spaces with the coarse wedge topology.

\smallskip
\noindent{\bf 3.11. Construction.}  Given a base $\mathcal{B}$ for
a non-Archimedean space $X$ such that $\mathcal{B}$ is a tree under
reverse containment, let $\langle \mathcal{T}, \le \rangle$ be the 
chain 
completion of $\mathcal{B}$.  For each $x \in X$ let $\mathcal{B}(x)$
be the branch of all $B \in \mathcal{B}$ such that $x \in B$.
Then the map $f\colon X \rightarrow \mathcal{T}$ that takes $x$
to the supremum of $\mathcal{B}(x)$ in $\mathcal{T}$ is easily seen
to be an embedding with respect to the coarse wedge topology.
\smallskip

It has long been known that every non-Archimedean space is realizable as 
some subset of the set of all branches of some tree, endowed with a 
natural topology analogous to the definition of the Stone space of a 
Boolean algebra. One simply takes a tree base $\mathcal{B}$ and proceeds 
as above; usually, the tree $\mathcal{T}$ is not explicitly mentioned, 
only the tree $\mathcal{B}$ and the set of its branches.

The analogy with the Stone duality goes in the opposite direction, too. 
Given a tree $S$, one can let the space $\mathcal{X}(S)$ the set of the 
branches of $S$.  The resulting space has a tree base $\mathcal{B}$ in 
natural association with $S$, with $s \in S$ corresponding to $B[s] = \{X 
\in \mathcal{X}(S) : s \in X\}$. See [21] for details.

Many properties of $S$ are naturally associated to properties of 
$\mathcal{X}(S)$.  For example, $\mathcal{X}(S)$ is an L-space if, and 
only if, $S$ is a Souslin tree [Definition 4.10 below].  One also has some 
carry-over in Construction 3.11, though one needs to be careful. If $X$ is 
an L-space, then every tree-base for $X$ is indeed a Souslin tree, but its 
completion $\mathcal{T}$ will not be an L-space if a finitary tree-base is 
chosen, since then $\mathcal{T}$ has uncountably many isolated points.  
On the other hand, if the base is chosen [as indeed it can be] so that 
every member, other than an isolated singleton, has infinitely many 
immediate successors, then $\mathcal{T}$ will be a compact L-space, as 
will be shown in a forthcoming paper.

An interesting class of non-Archimedean spaces is provided by
trees [and not their branch spaces!] in which each member has
at most countably many immediate successors, with the fine
wedge topology. For each  $t \in T$ with infinitely many immediate
successors, let $\langle t_n : n \in \w \rangle$ list them, and let
$B_t^n = V_t \setminus (V_{t_1} \cap \dots \cap V_{t_n})$.
Each $t \in T$ with finitely many immediate successors is isolated,
so that 
$$
\{\{t\} : t \text{ is isolated }\} \cup \{B_t^n : n \in \w, \
t \text{ has infinitely many immediate successors}\} 
$$ 
is a tree base for $T$ with the fine wedge topology.  One 
consequence of all this is something that foreshadows a theme
of the next section:

\begin{myproclaim}{3.12. Theorem} Every tree in which each element has
at most countably many immediate successors is suborderable in
the fine wedge topology.
\end{myproclaim}

Indeed, every non-Archimedean space is suborderable. One can also show 
that every full $\w$-ary tree of limit order height is orderable in the 
fine wedge topology. C. Aull [1] took advantage of this to produce a 
hereditarily paracompact space with a point-countable base with no 
$\sigma$-point-finite base, using the full $\w$-ary tree of height $\w_1$.  
This idea generalizes to all cardinals in a straightforward way.  
Incidentally, it is not hard to show that the Michael line is homeomorphic 
to the full $\w$-ary tree of height $\w + 1$ in the fine wedge topology, 
with the points at level $\w$ corresponding to the irrationals. Details 
will appear in a forthcoming paper.

{\it From now on, ``tree'' will always mean, ``tree in which every
nontrivial pseudo-supremum is a supremum.''}

\section*{A Short Survey on (Mostly) the Interval Topology}

The interval topology has received the lion's share of attention among 
set-theoretic topologists as far as topological properties are concerned.  
Part of the explanation for this is twofold: on the one hand, most of the 
topologies in Section 2 are not Hausdorff except in trivial cases; and on 
the other hand, the remaining ones (except for the fantail topology, which 
coincides with the interval topology on finitary trees) have such strong 
topological properties that there is far less room for variation than with 
the interval topology.  The following concept highlights this difference:

\smallskip
\noindent{\bf 4.1. Definition.}  A space $X$ is {\it monotone normal}
(or: {\it monotonically normal \/}) if to each pair $\langle G, x \rangle$
where $G$ is an open set and $x \in X$, it is possible 
to assign an open set $G_x$ such that $x \in G_x \subset G$ so that
$G_x \cap H_y \ne \emptyset$ implies either $x \in H$ or
$y \in G$.
\smallskip

[The foregoing is actually a characterization due to C.~R.~Borges [4] 
which is very well adapted to our purposes.  The usual definition 
motivates the name ``monotone normal'' much better.]

Monotone normality is a very strong property.  It is hereditary, and it 
implies both collectionwise normality and countable paracompactness.  So 
the following theorem tells us that trees are `very nicely behaved' in 
three of the first four topologies:

\begin{myproclaim}{\bf 4.2. Theorem}  Every tree is monotonically normal in
the coarse wedge, fine wedge, and chevron topologies.
\end{myproclaim}

\begin{proof}[Outline of Proof]
For the chevron topology, given $t \in G$,
let $G_t = \{t\}$ if $t$ is isolated, and otherwise let
 $G_t =  C[s, t]$ for the minimal $s$ such that $C[s,t] \subset G$.
For the fine wedge topology, $ G_t$ can be defined by removing from $V_t$
all of the finitely many $V_x$ that are not subsets of $G$ among
the immediate successors $x$ of $t$.  For the coarse wedge topology,
put the choices for the two other topologies together. 
\end{proof}

The well-known Rudin-Balogh characterization of [hereditary] 
paracompactness in monotonically normal spaces translates very simply to 
trees in these three topologies: a tree is paracompact iff it does not 
have a closed copy of an uncountable regular cardinal and hereditarily 
paracompact iff it has no copies of stationary subsets of uncountable 
regular cardinals. Since there are no such subspaces in the fine wedge 
topology at all, we have:

\begin{myproclaim}{\bf 4.3. Corollary} Every tree is hereditarily
paracompact in the fine wedge topology.
\mbox{} \hfill $\square$
\end{myproclaim}

The situation is completely different for the interval topology, where 
monotone normality imposes a very strong structure on the tree: it is 
equivalent to the tree being a topological direct sum of copies of ordinal 
spaces (Theorem 4.7 below).  This rules out such well-known examples as 
Aronszajn trees and the Cantor tree.

{\it For the rest of this article, all topological statements concerning 
trees will refer to the interval topology.}

Two other characterizations of monotone normal trees are given in the 
following two definitions.

\smallskip
\noindent{\bf 4.4. Definition.}  Let $\Lambda$ denote the class of limit
ordinals. A tree $T$ has {\it Property $\delta$} 
if there exists a function $f\colon T \restriction \Lambda \rightarrow T$ 
such that $f(t) < t$ for all $t \in T \restriction \Lambda$, and such 
that if 
$[f(s), s]$
meets $[f(t), t]$ then $s$ and $t$ are comparable.  

\smallskip

\noindent{\bf 4.5. Definition.} A {\it neighbornet} in a space
$X$ is a function $U\colon X \rightarrow \mathcal{P}(x)$ such that 
$U(x)$ is a neighborhood of $x$ for all $x \in X$.
 A neighbornet $V$ {\it refines} $U$ if $V(x) \subset U(x)$
for all $x \in X$. A space $X$ is {\it halvable} if each
neighbornet $U$ of $X$ has a neighbornet $W$ refining it such that
if $W(x) \cap W(y) \ne \emptyset$ then either $x \in U(y)$ or
$y \in U(x)$.

\smallskip

\noindent{\bf 4.6. Definition.}  A subset $S$ of a tree $T$
is {\it convex} if $[s_1, s_2] \subset S$ whenever $s_1$ and $s_2$ 
are elements of $S$ satisfying $s_1 < s_2$.
\smallskip

\begin{myproclaim}{4.7.  Theorem} {\rm [22]} 
Let $T$ be a tree.  The following are equivalent.
\begin{enumerate}
\item $T$ is monotonically normal.
\item $T$ is halvable.
\item $T$ has Property $\d$.
\item $T$ is the topological direct sum of subspaces, each homeomorphic
to an ordinal and each convex in $T$.
\item $T$ is orderable.
\item The neighborhoods of the diagonal in $T^2$ constitute a uniformity.
\hfill $\square$
\end{enumerate}
\end{myproclaim}

Paracompactness is even more restrictive.  Locally compact, paracompact, 
zero-dimensional spaces are the topological direct sum of compact clopen 
subspaces.  Hence, a tree is easily seen to be paracompact if, and only 
if, it is a topological direct sum of compact spaces, each homeomorphic to 
an ordinal.  Of course, this implies they are monotone normal.  Also, it 
is easy to see:

\begin{myproclaim}{4.8. Theorem}  The following are equivalent for a tree $T$:
\begin{enumerate}
\item $T$ is hereditarily paracompact.
\item $T$ is paracompact and has no uncountable branches.
\item $T$ is the topological direct sum of countable, compact spaces each 
homeomorphic to an ordinal.
\item  $T$ is metrizable.  
\hfill $\square$
\end{enumerate}
\end{myproclaim}

And so, most of the topological action here has to do with concepts weaker 
than monotone normality.  Many of these properties have been studied by 
set-theoretic topologists, but usually only in connection with what are 
rather cryptically called ``{\it $\w_1$-trees}''.  These are trees of 
height $\w_1$ in which every level is countable.  Usually even more 
conditions are imposed, such as the conditions that every element has 
successors at all levels of the tree and every element has at least two 
immediate successors; trees satisfying these latter two properties are 
often called {\it normalized}.

Strangely enough, the proofs of most of the general theorems in the 
literature about topological properties on $\w_1$-trees go through almost 
verbatim for arbitrary (Hausdorff, by the conventions of these last two 
sections) trees.  One of the rare exceptions is Theorem 4.7 above, where 
the proof that (3) implies (4) in [16] really does not adapt readily to 
the general case.  In some of the theorems below, however, I will not even 
add ``in effect'' when attributing them to various authors, so close is 
the published proof to one for trees in general.  This applies to the 
following theorem, which introduces an important motif: many familiar 
topological properties can be reduced to the case where all or all but one 
of the initial ingredients is an antichain.

\begin{myproclaim}{4.9. Theorem} {\rm Fleisner, [10]}
Let $T$ be a tree.  The following are equivalent.
\begin{enumerate}
\item $T$ is normal.
\item Given a closed set $F$ and an antichain $A$ disjoint from $F$,
there are disjoint open sets $G$ and $H$ such that $A\subset G$
and $B \subset H$. 
\hfill $\square$
\end{enumerate}
\end{myproclaim}

Some of the most important classes of trees have definitions involving 
antichains.

\smallskip
\noindent{\bf 4.10. Definition.}  A tree is {\it special} if 
it is a countable union of antichains.  A tree is {\it Souslin}
if it is uncountable while every chain and antichain is countable.
A tree is {\it Aronszajn} if it is uncountable while every chain
is countable and every level $T(\a)$ is countable.
\smallskip

One of the most useful and obvious topological facts about trees is that 
every antichain is a closed discrete subspace.  A closely related result 
is:

\begin{myproclaim}{4.11. Theorem} Let $X$ be a subset of a tree $T$.  The
following are equivalent:
\begin{itemize}
\item[(i)] $X$ is a countable union of antichains.
\item[(ii)] $X$ is $\sigma$-discrete, i.e., it is a countable union of
closed discrete subspaces.
\end{itemize}
\end{myproclaim}

\begin{proof}[Proof that (ii) implies (i):]
It is clearly
enough to show that every closed discrete subspace is 
the countable union of antichains.  So let $D$ be closed
discrete, let $D_0$ be the set of minimal members of $D$,
and with $D_n$ defined, let $D_{n+1}$ be the set of minimal
members of $D \setminus (D_0 \cup \dots \cup D_n)$.
Clearly each $D_n$ is an antichain of $T$.  
If there were a point $d$ in $D$ but not any of the
$D_n$, then for each $n \in \w$ there would be a point 
$d_n \in D_n$ such that $d_n < d$, and any point in the
pseudo-supremum of the $d_n$ would be in their closure,
violating the claim that $D$ is closed discrete.
\end{proof}

Thus, in particular, every special tree is a countable union of closed 
discrete subspaces.  This clearly implies every chain is countable, and 
hence also easily implies that each special tree is developable.  This was 
shown by F.~Burton Jones, who gave special Aronszajn trees the name ``tin 
can spaces,'' investigating them over a period of many years as candidates 
for a nonmetrizable normal Moore space, along with the related ``Jones 
road spaces'' formed from them in the way described near the end of 
Section 1. His judgment was partially vindicated when W. Fleissner showed 
[9] that these spaces are normal under MA + $\neg$CH. However, Devlin and 
Shelah [6] showed that no special Aronszajn tree is normal under 
$2^{\aleph_0} < 2^{\aleph_1}$. Ironically enough, this was the same axiom 
that Jones used back in 1937 to show the consistency of every separable 
normal Moore space being metrizable.  Thus, in particular, the situation 
as regards ``$\w_1$-Cantor trees'' and special Aronszajn trees is exactly 
parallel: the trees obtained by removing all except exactly $\w_1$ points 
from the top level of the Cantor tree are nonmetrizable Moore spaces, as 
are special Aronszajn trees; MA($\w_1$) implies both classes of trees are 
normal; and $2^{\aleph_0} < 2^{\aleph_1}$ implies both classes are not 
normal.

Special trees have another property, which is often given the name 
``$\mathbb{Q}$-embeddability''; but the map involved is almost never 
a topological embedding, nor is it usually one-to-one.  So the following 
terminology is adopted here:

\smallskip
\noindent{\bf 4.12. Definition.}  Let $\langle L, \le_L \rangle$ 
be a linearly ordered set, and let $\langle P, \le_P \rangle$ 
be a tree.  A function $f\colon P \rightarrow L$ is called an
{\it $L$-labeling} if it is strictly order preserving: that is,
$p <_P q$ implies $f(p) <_L f(q) $.  A poset is {\it $L$-special}
if it admits an $L$-labeling. 

\begin{myproclaim}{4.13. Theorem}  Let $T$ be a tree.  The following are equivalent:
\begin{itemize}
\item[(i)] $T$ is special.
\item[(ii)] $T$ is $\mathbb{Q}$-special.
\item[(iii)] $T$ is $\sigma$-discrete in the interval topology.
\item[(iv)] $T$ is developable in the interval topology.
\item[(v)] $T$ is subparacompact in the interval topology, and is of 
height $\le \w_1$.  
\hfill $\square$
\end{itemize}
\end{myproclaim}

The proof that (i) is equivalent to (ii) is well known (cf.~9.1 of [28].  
The equivalence of (iii) through (v) was demonstrated, in effect, by K.-P. 
Hart [15] although the class of trees explicitly mentioned was more 
restrictive.

\smallskip
\noindent{\bf 4.14. Definition.} Let $\mathcal{A}$ be a collection 
of disjoint nonempty sets.  A family $\mathcal{U}$ of sets {\it expands} 
$\mathcal{A}$ if for each $A \in \mathcal{A}$ there exists $U_A \in 
\mathcal{U}$ such that $A \subset U_A$ and $B \cap U_A = \emptyset$ if $B 
\ne A$. In case where $\mathcal{A}$ consists of singletons, we also say
$\mathcal{U}$ {\it expands} $\bigcup \mathcal{A}$.

\smallskip

\noindent{\bf 4.15. Definition.} A space $X$ is [{\it strongly}] {\it collectionwise
Hausdorff } (often abbreviated {\it [s]cwH}) if every closed discrete
subspace expands to a disjoint [{\it resp.} discrete] collection
of open sets.  A space $X$ is {\it collectionwise
normal} (often abbreviated {\it cwn}) if every discrete collection
of closed sets expands to a disjoint (equivalently, discrete) collection
of open sets.

\begin{myproclaim}{4.16. Theorem} {\rm (M.~Hanazawa [14])} 
Let $S$ be a subspace of a tree.  The following
are equivalent.
\begin{enumerate}
\item $S$ is collectionwise Hausdorff (cwH).
\item Every antichain of $S$ expands to a disjoint collection of open sets.
\item $S$ is hereditarily cwH.
\hfill $\square$
\end{enumerate}
\end{myproclaim} 

\begin{myproclaim}{4.17. Theorem} {\rm  (K.~P.~Hart [16, proof of 2.1],  
in effect)}  Let $S$ be a subspace of a tree.  The following are 
equivalent:
\begin{enumerate}
\item S is normal and cwH.
\item S is strongly cwH.
\item S is hereditarily collectionwise normal.  
\hfill $\square$
\end{enumerate}
\end{myproclaim}

\begin{myproclaim}{4.18. Corollary}  Every Souslin tree is hereditarily
collectionwise normal.
\end{myproclaim}

\begin{proof}
Every antichain $A$ is countable and hence is a subset
of some clopen initial segment $T \restriction (\a + 1)$, 
which is second
countable and thus metrizable. Therefore, $A$ can be expanded to
a discrete collection of open subsets of $T \restriction (\a + 1)$, 
and
hence of $T$.
\end{proof}

\begin{myproclaim}{4.19. Corollary} The existence of a normal tree that is
not collectionwise normal is ZFC-independent.
\end{myproclaim}

\begin{proof}
If MA + $\neg$CH, one can either use a special 
Aronszajn tree (which is not cwH by the Pressing-down Lemma)
or an $\w_1$-Cantor tree, as remarked early in Section 4, to give such 
a tree.  On the other hand $V=L$ implies every locally compact
normal space is cwH [29], and so Theorem 4.17
implies it is (hereditarily) cwn.
\end{proof}

For our next few results, recall that a space is said to be {\it countably 
paracompact} [resp.\ {\it countably metacompact\/}] if every countable 
open cover has a locally finite [resp. point-finite] open refinement.

\begin{myproclaim}{4.20. Theorem} {\rm (Nyikos [24])} Let $T$
be a tree.  The following are equivalent:
\begin{enumerate}
\item $T$ is countably paracompact {\rm [resp. {\it countably 
metacompact \/}]}.
\item Any countable partition $\{A_n : n \in \w\}$ of any antichain
of $T$ expands to a locally finite {\rm [resp. {\it point-finite}]}
collection of open sets.
\hfill $\square$
\end{enumerate}
\end{myproclaim}

\begin{myproclaim}{4.21. Corollary} 
Every normal tree is countably paracompact. (``There are no Dowker 
trees.'') 
\end{myproclaim}

\begin{proof}
Every countable discrete collection of closed sets
in a normal space expands to a discrete collection of open sets.
\end{proof}

\begin{myproclaim}{4.22. Corollary} 
Every cwH tree is countably metacompact.
\end{myproclaim}

\begin{proof}
Given $\{A_n : n \in \w\}$ as in 4.20, use cwH to
expand the antichain that is their union to a disjoint collection
of open sets, and let $U_n$ be the union of the ones that meet
$A_n$.
\end{proof}

\begin{myproclaim}{4.23. Theorem} The existence of a countably paracompact
tree that is not cwH is ZFC-independent. \end{myproclaim}

\begin{proof}
W.~S.~Watson showed that under $V = L$,
every locally compact, countably paracompact space is cwH [30].  On
the other hand, a $\w_1$-Cantor tree is not cwH, but is normal
under MA + $\neg$CH, and hence countably paracompact under the
same axiom because it is normal and Moore.
\end{proof}

The following three questions are related to the last three results.  
Note the contrast in the phrasing as to set-theoretic status.

\begin{myproclaim}{4.24. Problems} Is there a ZFC example of a cwH
tree that is {\bf (a)} not countably paracompact or 
{\bf (b)} not normal or at least {\bf (c)} not monotone normal?
\end{myproclaim}

{\it Caution.} 
K.~P.~Hart  credits S.~Todor\v cevi\'c with having even shown, assuming 
``an at least inaccessible cardinal'', that  it is consistent for every 
cwH tree to be orderable [16].  If this had been correct as stated, these 
problems would be solved modulo inaccessibles, but ``tree'' referred to 
$\w_1$-trees only.

If there is a Souslin tree, there is a cwH tree which is not countably 
paracompact and hence (by 4.21) not normal.  Details will appear in [24].  
Earlier, Devlin and Shelah [7] used the stronger axiom $\diamondsuit^+$, a 
consequence of V = L, to construct a cwH non-normal tree which is 
$\mathbb{R}$-special, hence not countably paracompact (see Corollary 
4.40 below).

\begin{myproclaim}{4.25. Problem} Is it true in ZFC that every 
countably paracompact cwH tree is (collectionwise) normal?
\end{myproclaim}

\begin{myproclaim}{4.26. Problems} Is it ZFC-equiconsistent that 
every countably paracompact tree is {\bf (a)} normal?
{\bf (b)} collectionwise normal?
\end{myproclaim}

The last question is phrased the way it is because of a gap in our 
consistency results.  On the one hand, a $\Delta$-set of real numbers that 
is not a $Q$-set gives a countably paracompact non-normal tree, and such 
sets of reals are consistent assuming just the consistency of ZFC [17]. On 
the other hand, the only known models in which every locally compact, 
countably paracompact space is {\it strongly} cwH require large cardinal 
axioms.  This also applies to ``first countable'' in place of ``locally 
compact,'' and we know of no axioms which give normality in countably 
paracompact trees without also giving cwH.  In fact, the following problem 
is of interest quite apart from its obvious applicability to Problem 4.26:

\begin{myproclaim}{4.27. Problem} Does V = L or some other ZFC-equiconsistent
axiom imply that every locally compact,
or every first countable, countably paracompact regular space
is strongly cwH?  \end{myproclaim}

Unlike countable paracompactness, countable metacompactness has generally 
been thought of as a very weak property.  However, the following suggests 
that its failure for trees is a fairly ordinary occurrence:

\begin{myproclaim}{4.28. Example} {\rm (Nyikos [24])} The full
binary tree of height $\w_1$ is not countably metacompact. \end{myproclaim}

There are even trees that are $\mathbb{R}$-special, yet not countably 
metacompact, such as the tree of all ascending sequences of rational 
numbers [24], designated $\sigma \mathbb{Q}$ in [28].

Being $\mathbb{R}$-special turns out to be a fairly strong ``generalized 
metric'' property for trees.  It is easily shown to imply 
quasi-metrizability, but is much stronger [23], and we also have:

\begin{myproclaim}{4.29. Theorem} {\rm (K.~P.~Hart, [15])}  Let
$T$ be a tree.  The following are equivalent.
\begin{enumerate}
\item $T$ is $\mathbb{R}$-special.
\item $T$ has a $G_\d$-diagonal.
\item The set of nonisolated points of $T$ is a $G_\d$.
\item The set of isolated points of $T$ is a countable union of antichains.
\hfill $\square$
\end{enumerate}
\end{myproclaim}

In the same article, Hart also showed the remarkable fact that every 
finitary $\mathbb{R}$-special tree is special.  So, for example, if each 
element of $T$ has $\le {\mathfrak c}$-many immediate successors, and we 
add a full binary tree of height $\w$ between each point and its immediate 
successors, then $T$ embeds as a closed subspace in the resulting tree, 
and if $T$ is special, so is the resulting tree.  On the other hand, if 
$T$ is $\mathbb{R}$-special but not special, the resulting tree will be 
quasi-metrizable but not special [23]. Much is still unknown about 
quasi-metrizable trees, including:

\smallskip
\noindent{\bf 4.30 Problem.} Is it consistent that every tree
without uncountable branches is quasi-metrizable?
\smallskip

This problem is worded the way it is because $\w_1$ embeds in
every tree with an uncountable branch, and is not quasi-metrizable;
and because a Souslin tree is not quasi-metrizable [23].
 
Condition (4) in Theorem 4.29 was an ingredient in the proof of:

\begin{myproclaim}{4.31. Theorem} {\rm ( Nyikos [24])} Let $T$ be
a tree.  The following are equivalent:
\begin{enumerate}
\item $T$ is perfect; that is, every closed subset of $T$
is a $G_\d$.
\item $T$ is $\mathbb{R}$-special, and every antichain is a $G_\d$.
\hfill $\square$
\end{enumerate}
\end{myproclaim}

The fact that (1) implies (2) was essentially shown by M.~Hanazawa [14]] 
who used it to help answer a question of K.~P.~Hart [15]: ``Is every 
$\w_1$-tree with a $G_\d$-diagonal perfect?'' The answer is affirmative 
under MA\,+$\,\neg$CH as Hart himself pointed out, but it is negative 
under the axiom $\diamondsuit^*$ with the help of which Hanazawa 
constructed a counterexample.  [{\it Caution.} The example, an Aronszajn 
tree, is claimed in [13] to be countably metacompact, but it is not.] The 
paper also showed the following for $\w_1$-trees:

\begin{myproclaim}{4.32. Corollary} 
Every collectionwise Hausdorff, $\mathbb{R}$-special tree is perfect.
\end{myproclaim}

\begin{proof}
Because $T$ is $\mathbb{R}$-special, its height is $\le \w_1$.
Hence, by cwH, every antichain expands to a disjoint family of countable
open sets, and hence is a $G_\d$.
\end{proof}

Under MA + $\neg$CH, we can weaken the hypothesis and strengthen the
conclusion:

\begin{myproclaim}{4.33. Theorem} {\rm (Nyikos [22])} If MA + $\neg$CH,
a tree is metrizable if, and only if, it is collectionwise Hausdorff
and has no uncountable chains. \end{myproclaim}

A related ZFC result is:

\begin{myproclaim}{4.34. Theorem} {\rm (Nyikos [22])} A tree is 
metrizable if, and only if, it is special and cwH. \end{myproclaim}

Another corollary of Theorem 4.31 is:

\begin{myproclaim}{4.35. Corollary} If V = L, or PMEA, then the following are equivalent
for a tree $T$:
\begin{enumerate}
\item $T$ is perfect.
\item $T$ is $\mathbb{R}$-special and countably metacompact.
\end{enumerate}
\end{myproclaim}

\begin{proof}
As is well known, every perfect space is countably
metacompact, so from 4.31 it follows that (1) implies (2).  Conversely,
if V = L, then every closed discrete subspace
in a locally countable, countably metacompact space is a $G_\d$ [20]; and
every $\mathbb{R}$-special tree is of height $\le \w_1$ and hence locally
countable.  Under PMEA, every closed discrete subspace in a first
countable, countably metacompact space is a $G_\d$ (D.~Burke,
[5]) and (2) similarly implies (1). 
\end{proof}  

I do not know whether the set-theoretic hypotheses in 4.35 can be dropped.  
More generally:

\begin{myproclaim}{4.36 Problem} Is every closed discrete subset of a countably
metacompact tree of height $\le \w_1$ a $G_\d$?
\end{myproclaim}

The following is a pleasing counterpoint to Theorem 4.31:

\begin{myproclaim}{4.37. Theorem} ({\rm Nyikos [24])}  Let $T$ be a tree.
The following are equivalent.
\begin{enumerate}
\item $T$ is perfectly normal.
\item Every closed subset of $T$ is a regular $G_\d$; that is,
it is of the form $\bigcap \{cl(U_n) : n \in \w\}$ where each
$U_n$ is open.
\item $T$ is $\mathbb{R}$-special and every antichain is a regular $G_\d$.
\hfill $\square$
\end{enumerate}
\end{myproclaim}

Perfect normality is a highly axiom-sensitive property where trees are 
concerned.  Under V = L, normal trees are collectionwise normal, and 
M.~Hanazawa [14] used the consequence $\diamondsuit^*$ of V = L to 
construct a perfectly normal Aronszajn tree which is, of course, not 
special.  Under MA + $\neg$CH, the cwH ones are all metrizable (Theorem 
4.32), but every Aronszajn tree is a nonmetrizable example, and special, 
as is every $\w_1$-Cantor tree.  Under $2^{\aleph_0} < 2^{\aleph_1}$ and 
EATS (``Every Aronszajn tree is special'') no Aronszajn tree can be normal 
(cf.~[6] or [27]), nor can any $\w_1$-Cantor tree, by the Jones Lemma.  
However, there is an axiom compatible with CH under which there is a 
perfectly normal non-cwH tree of height $\w + 1$ [22]. Finally, if 
strongly compactly many random reals are added to a model of MA + 
$\neg$CH, Theorem 4.32 still holds in the forcing extension, but PMEA also 
holds and so every normal, first countable tree is cwn, hence every 
perfectly normal tree is metrizable.  However, the following may still be 
open:

\begin{myproclaim}{4.38. Problem} Is it ZFC-equiconsistent for every perfectly
normal tree to be metrizable? \end{myproclaim}

Weakening normality slightly to countable paracompactness, we have:

\begin{myproclaim}{4.39. Theorem} 
({\rm Nyikos [24])}  Every $\mathbb{R}$-special, cwH, countably 
paracompact tree is (collectionwise) normal.  
\hfill $\square$ 
\end{myproclaim}

\begin{myproclaim}{4.40. Corollary} If V = L, every $\mathbb{R}$-special, 
countably paracompact tree is collectionwise normal.
\end{myproclaim}

\begin{proof}
By 4.39 and the proof of 4.24. 
\end{proof}

Thus Problems 4.25 and 4.26 have affirmative answers for 
$\mathbb{R}$-special trees.  

\begin{myproclaim}{4.41. Problems} Is there a ZFC example of a tree which
is not special but is {\bf (a)} perfect or {\bf (b)} countably
metacompact and has no uncountable branches? \end{myproclaim}

The search for ZFC examples is made difficult by the fact that $\sigma 
\mathbb{Q}$, which is the simplest ZFC example of an $\mathbb{R}$-special, 
non-special tree of which I am aware, is not countably metacompact.  
Consistent examples have long been known, like the Devlin-Shelah 
$\diamondsuit^*$ example mentioned earlier: a cwH (hence non-special, by 
the Pressing-Down Lemma) Aronszajn tree which is not normal but is 
$\mathbb{R}$-special, hence perfect.  In the same paper [7], they 
defined:

\smallskip
\noindent{\bf 4.42. Definition.} An Aronszajn tree is {\it almost
Souslin} if every antichain meets a nonstationary set of levels. 
\smallskip

They showed that an Aronszajn tree is cwH iff it is almost Souslin; this 
is an easy application of Theorem 4.11 and the Pressing-down Lemma.

Back in the 1980's, Hanazawa did some extensive cataloguing of how 
Aronszajn trees behave under V = L.  He constructed or listed examples 
exhibiting all combinations of the properties considered in this section 
and not ruled out by the results mentioned or proven here --- except one: 
we still do not know whether there is a countably paracompact, non-normal 
Aronszajn tree under V = L. He also has catalogued their behavior with 
respect to some properties not mentioned here, cf.~[12].

Finally, we return to the property with which we began this section.

\begin{myproclaim}{4.43. Problems} Is it consisent that every {\bf (a)} normal
or {\bf (b)} every collectionwise normal tree is monotonically normal?
\end{myproclaim}

Recall also Problem 4.24 (c), which can be stated negatively: is it 
consistent that every cwH tree is monotonically normal? Here is a partial 
result:

\begin{myproclaim}{4.44 Theorem} {\rm(Nyikos [22])} If MA + $\neg$wKH, 
then every cwH tree of height $< \w_2$ is monotonically normal. 
\end{myproclaim}

[Compare Theorem 4.33.] Here ``wKH'' refers to the existence of weak 
Kurepa trees. Once informally called ``Canadian trees,'' these are trees 
of height and cardinality $\w_1$ that have more than $\aleph_1$ 
uncountable branches.  Of course, the full binary tree of height $\w_1$ is 
a weak Kurepa tree under the continuum hypothesis.  So the axiom in 4.44 
negates CH, and it is known to imply that there are inaccessible 
cardinals.
 
It would be interesting to see whether large cardinals are really needed 
in 4.44.  Its proof goes through just on the assumption of ``every weak 
Kurepa tree has a special Aronszajn subtree,'' whose consistency is 
apparently not known to depend on large cardinal axioms.

One of the ingredients in the proof of 4.44 is of independent interest.  
Call a tree {\it $\sigma$-orderable} if it is the countable union of 
closed, orderable subtrees.

\begin{myproclaim}{Theorem 4.45} {\rm (Nyikos [22])}
A tree is orderable if, and only if, it is $\sigma$-orderable and cwH. 
\hfill $\square$
\end{myproclaim}

Incidentally, $\sigma$-orderability is easily seen to be equivalent to a 
property to which J.~E.~Baumgartner tried to transfer the term ``special'' 
[3].  His usage does not, however, seem to have caught on, and 
``$\sigma$-orderable'' seems to be a good a name as any for Baumgartner's 
property.

\section*{Acknowledgements}

I am grateful to my graduate students Chunliang Pan and Akira Iwasa for 
their help with this manuscript. Mr.~Pan produced the Postscript files for 
Figures 1 and 2, and Mr.~Iwasa caught several errors in earlier drafts of 
Section 2.  I would also like to thank M.~Hanazawa and Klaas~Pieter~Hart 
for sending me reprints of their articles, which not only helped give body 
to Section 4 but also spurred me on to much of the research that led to 
this paper.


\begin{thebibliography}{30}

\bibitem{1}
C. E. Aull, 
{\it Topological spaces with a $\sigma$-point-finite base,} AMS 
Proceedings {\bf 29:}411--416 (1971).

\bibitem{2}
B. Balcar, J. Pelant, and P. Simon, {\it The space of ultrafilters on N 
covered by nowhere dense sets,}  Fund. Math. {\bf 110}:11-24 (1980).

\bibitem{3}
J. E. Baumgartner, {\it Applications of the 
Proper Forcing Axiom}, in: ``Handbook 
of Set-Theoretic Topology,'' K. Kunen and J. Vaughan, eds., North-Holland, 
Amsterdam, 913--959 (1984).

\bibitem{4}
C. R. Borges, {\it Four generalizations of 
stratifiable spaces,} in: ``General 
Topology and its  Applications to Modern Analysis and Algebra III,'' Akademia,
Praha,  
 73--77 (1971). 

\bibitem{5}
D. K. Burke, {\it PMEA and first countable, countably paracompact spaces,}
AMS Proceedings {\bf 92}:455-460 (1984).

\bibitem{6}
K. J. Devlin and S. Shelah,  {\it A note on the normal Moore space 
conjecture, } Canad. J. Math. {\bf 31:}241--251 (1979).

\bibitem{7}
K. J. Devlin and S. Shelah, {\it Suslin properties and tree topologies,}  
Proc. London Math. Soc. {\bf 39:} 237--252 (1979).

\bibitem{8}
J. Diestel,  ``Sequences and series in Banach spaces,''
Springer-Verlag (1984).

\bibitem{9}
W. G. Fleissner, {\it When is Jones' space normal? }
AMS Proceedings {\bf 50:}375--378 
(1975).

\bibitem{10}
W. G. Fleissner, {\it Remarks on Suslin properties and 
tree topologies,} AMS Proceedings {\bf 80:}320--326 (1980).

\bibitem{11}
G. Gruenhage,
{\it On a Corson compact space of Todor\v cevi\'c,} Fund. Math.
{\bf 126}:261--268 (1986).

\bibitem{12}
M. Hanazawa,
{\it Various kinds of Aronszajn trees with no subtree of a different
kind,} in: Lecture Notes in Mathematics  \#891, Springer-Verlag, 
1--21 (1981).

\bibitem{13}
M. Hanazawa, {\it Countable metacompactness and tree topologies,} J. Math. 
Soc. Japan {\bf 35}: 59--70 (1983).

\bibitem{14}
M. Hanazawa,
{\it Note on countable paracompactness of collectionwise Hausdorff
tree topologies,} Saitama Math. J. {\bf 2}:7--20 (1984).

\bibitem{15}
K. P. Hart,
{\it Characterizations of $\mathbb{R}$-embeddable and developable
 $\omega_1$-trees,} Indag. Math. {\bf 44}:277--283 (1982).

\bibitem{16}
K. P. Hart,
 {\it More remarks on {Souslin} properties and tree topologies,}
 Top. Appl. {\bf 15}: 151--158 (1983).

\bibitem{17}
R. W. Knight, {\it $\Delta$-sets,} AMS Transactions {\bf 339}:45--60.

\bibitem{18}
J. D. Lawson, {\it The versatile continuous order,} in: ``Mathematical
Foundations of Programming Language Semantics,'' ed. by
M. Main, A. Melton, M. Mislove, and D. Schmidt, Lecture 
Notes in Computer Science \#298, Springer-Verlag, 134--160 (1988).

\bibitem{19}
J. W. Morgan, {\it $\Lambda$-trees and their applications}, 
AMS Bulletin {\bf 26}: 87--112 (1992).

\bibitem{20}
P. J. Nyikos, {\it Countably metacompact, locally countable spaces in the 
constructible universe,} Coll. Math. Soc. J\'anos Bolyai {\bf 
55:}409--424 (1989).

\bibitem{21}
P. J. Nyikos,
{\it On a paper of Alexandroff and Uryshohn on certain non-Archimedean 
spaces,} Top. Appl. (to appear).

\bibitem{22}
P. J. Nyikos,
{\it Metrizability and orderability of trees} (in preparation).

\bibitem{23}
P. J. Nyikos,
{\it Quasi-metrizability in trees and ordered spaces} (in preparation).

\bibitem{24}
P. J. Nyikos,
{\it Countable paracompactness, countable metacompactness, and 
continuous layering of trees} (in preparation).

\bibitem{25}
D. Scott,
{\it Continuous lattices,} in: ``Toposes, Algebraic Geometry, and
Logic,''  Lecture Notes in Mathematics \#274, Springer-Verlag (1972).

\bibitem{26}
L. A. Steen and J. A. Seebach,  ``Counterexamples in Topology'' (Second 
Edition), Springer-Verlag (1978).

\bibitem{27}
A. D. Taylor,
{\it Diamond principles, ideals, and the normal Moore space problem,}
 Canad.  J. Math. {\bf 33}:282-296 (1981).

\bibitem{28}
S. Todor\v cevi\'c,
{\it Trees and linearly ordered sets,} in: ``Handbook of Set-Theoretic 
Topology,'' K. Kunen and J. Vaughan, eds. (North-Holland, Amsterdam),
235--293 (1984).

\bibitem{29}
W. S. Watson, {\it Locally compact  normal spaces in the constructible 
universe}, Canad. J. Math. {\bf 34}:1091--1096 (1982).

\bibitem{30}
W. S. Watson, {\it Separation in counatbly paracompact spaces,} AMS 
Transactions {\bf 290}:831--842 (1985).
\end{thebibliography}
\end{document}